\documentclass[11pt]{article}
\usepackage[T2A]{fontenc}
\usepackage[cp1251]{inputenc}
\usepackage{amsfonts}
\usepackage{eufrak}
\usepackage{amssymb}
\usepackage{amsmath}


\begin{document}

\centerline {\bf Independent random variables on Abelian}
\centerline{\bf groups with independent the sum and difference}

\bigskip

\centerline {\bf G.M. Feldman}

\bigskip

\begin{abstract}  Let $X$ be   a second
countable locally compact Abelian group. Let $\xi_1, \ \xi_2$ be
independent  random variables with values in the group $X$ and
distributions  $\mu_1, \ \mu_2$ such that the sum $\xi_1+\xi_2$
and the difference $\xi_1-\xi_2$ are independent. Assuming that
the connected component of zero of the group $X$ contains a finite
number elements of order 2 we describe the possible distributions
$\mu_k$.
\end{abstract}

\bigskip

{\it Key words.} Locally compact Abelian group,  Kac--Bernstein
theorem, Gaussian measure.

\bigskip

{\bf 1. Introduction.}
The classical Kac--Bernstein theorem states:

{\bf Theorem A.} {\it Let $\xi_1, \ \xi_2$ be independent random
variables. If the sum $\xi_1+\xi_2$ and the difference
$\xi_1-\xi_2$ are independent, then the random variables $\xi_k$ are Gaussian.}

Much research has been devoted to
group analogues of this  theorem
  (see \cite{Ru1} -- \cite{Fe6}).
In the present article we study the following question: Let
$\xi_1, \ \xi_2$ be independent random variables with values in a
locally compact Abelian group. Assume that the sum $\xi_1+\xi_2$
and the difference $\xi_1-\xi_2$ are independent. What one can say
about distributions of the random variables  $\xi_k$.
Before we pass to the proof of  the main result recall some necessary notation and
definitions.

Let $X$ be a  second countable locally compact Abelian group.
Denote by $Y=X^\ast$ the  character group of the group $X$, and by
$(x,y)$ the value of a character $y \in Y$ at a point $x \in
X$. Denote by $\mathbf{T}=\{z\in\mathbf{C}:\ |z|=1\}$ the circle
group (the one-dimensional torus).  Denote by ${M^1}(X)$ the
convolution semigroup of probability distributions on $X$. The characteristic function
of a distribution
$\mu \in {M^1}(X)$  we define  by the formula
 $$\widehat
\mu(y) = \int_X (x, y) d\mu(x).
$$
  A distribution $\gamma \in {M^1}(X)$
is called Gaussian (\cite{PaRaVa}), if its characteristic function
is represented in the form
\begin{equation}
\label{1}  \widehat\gamma(y)= (x,y)\exp\{-\varphi(y)\}, \quad y \in
Y,
\end{equation}
where $x \in X$ and $\varphi(y)$ is a continuous nonnegative
function on $Y$ satisfying the equation
\begin{equation}
\label{2} \varphi(u+v)+\varphi(u-v)=2[\varphi(u)+ \varphi(v)],
\quad u, \ v \in Y.
\end{equation}
Denote by $\Gamma(X)$ the set of Gaussian distributions on the group $X$.
A Gaussian distribution  $\gamma$ is called symmetric if
 in (\ref{1}) $x=0$. Denote by $\Gamma^s(X)$ the set of symmetric Gaussian distributions on $X$.
  Denote by $I(X)$ the set of idempotent distributions on the group $X$, i.e.
the set of shifts of
 Haar distributions $m_K$ of compact subgroups $K$ of the group $X$.

We will formulate now the following problem.

{\bf Problem 1}. Let $X$ be a  second countable locally compact Abelian group.
Let $\xi_1, \ \xi_2$ be  independent  random
variables with values in   $X$ and distributions  $\mu_1, \
\mu_2,$. Assume that the sum $\xi_1+\xi_2$ and the difference
$\xi_1-\xi_2$ are independent. Describe the possible distributions
$\mu_k$.

A.L. Rukhin in  \cite{Ru1} and  \cite{Ru2}  received some
sufficient conditions for the group $X$ in order that
distributions  $\mu_1$ and $\mu_2$ were represented as
convolutions of Gaussian and idempotent distributions. The
complete descriptions of such groups has been obtained by the
author.

{\bf Theorem B}  (\cite{Fe1}), see also  (\cite[\S 7]{Fe4}). {\it
Let $X$ be a   second countable locally compact Abelian group.
Assume that the connected component of zero of  the group $X$
 contains no elements of order $2$. Let   $\xi_1, \xi_2$ be
 independent random variables
with values in  $X$ and  distributions $\mu_1, \mu_2$ such that
the sum $\xi_1+\xi_2$ and the difference $\xi_1-\xi_2$ are
independent. Then $\mu_k\in \Gamma(X)*I(X), \ k=1, 2,$ and
$\mu_1=\mu_2*E_x, \ x\in X.$

 If the connected component of zero of a group $X$
 contains elements of order $2$, then
 there exist   independent random variables  $\xi_1, \xi_2$
 with values in  $X$ and  distributions $\lambda_1,
\lambda_2$   such that the sum $\xi_1+\xi_2$ and the difference
$\xi_1-\xi_2$ are independent, but
 $\lambda_k\notin \Gamma(X)*I(X), \ k=1, 2$.}

 We note that if a
distribution $\mu \in \Gamma(X)*I(X),$ then $\mu$ is invariant
with respect to a compact subgroup $K$, and $\mu$ induces a
Gaussian distribution on the factor-group $X/K$  under the natural
homomorphism $X \mapsto X/K$.

Problem 1 was solved in  \cite{BaESta} for the group
$X=\mathbf{T}$, and was solved in \cite{My} for the group
$X=\mathbf{R}\times \mathbf{T}$ and
$\text{\boldmath $a$}$-adic solenoids
 $\Sigma_{\text{\boldmath $a$}}$. Taking into account Theorem B solution of
 Problem 1 is reduced to the description of distributions $\mu_1$ and
  $\mu_2$ for groups $X$ which contain elements of order 2.
  In the present article  we   solve Problem 1 for groups
  $X$ satisfying the following condition: the connected component of zero of the group  $X$
 contains a finite number of elements of order 2. Note that if this condition holds,
 under additional assumption that the random variables $\xi_1$ and $\xi_2$ are identically
distributed Problem 1 was solved in \cite{Fe6}.

We need some results about the structure of locally compact
Abelian groups and the duality theory (see \cite[Ch. 6]{HeRo}).
  If $G$ is a closed subgroup of $X$,  then denote by
   $A(Y, G) = \{y \in Y: (x, y) = 1$ for all $x \in G \}$ its annihilator.
   The factor-group $Y/A(Y,G)$ is
topologically isomorphic to the character group of the group  $G$.
Put ${X_{(n)}=\{x\in X: nx=0 \}}, \ {X^{(n)}=\{x\in X:
\overline{}x=n\widetilde x, \ \widetilde x\in X\}}$. A group $X$ is said to
be Corwin group if $X^{(2)}=X$. Denote by
$\overline{Y^{(2)}}$ the closure of the subgroup ${Y^{(2)}}$.
Consider the subgroup $X_{(2)}$. It is
obvious that $A(Y, X_{(2)})=\overline{Y^{(2)}}$. Assume that
$X_{(2)}$ is a finite subgroup, and let $n$ be the number of its
elements $|X_{(2)}|=n$. Then $X_{(2)}\cong(X_{(2)})^*\cong Y/\overline{Y^{(2)}}$. Let $Y=\bigcup\limits_{j=0}^{n-1}
(y_j+\overline{Y^{(2)}})$ be  a
decomposition of the group $Y$ with respect to the subgroup $\overline{Y^{(2)}}$. If $G$ is a subgroup of
$X_{(2)}$, then its annihilator $H=A(Y, G)$, as is easily seen, is
of the form $H=\bigcup\limits_{j=0}^{l-1}
(y_j+\overline{Y^{(2)}})$. Denote by $c_X$ the connected component of zero
of the group
 $X$. Denote by $\mathbf{Z}$ the group of integers, and by
  $\mathbf{Z}(m)= \{0, 1,\dots, m-1\}$
the group of residue classes modulo $m$.

Let $f(y)$ be an arbitrary function on the group $Y$, and  let $h
\in Y$. Denote by $\Delta_h$ the finite difference operator
$$
\Delta_h f(y) = f(y+h) - f(y), \quad y \in Y.
$$

 Let $K$ be a compact subgroup of the group $X$. Then the characteristic function
 of the Haar distribution $m_K$ is of the form
 $$ \widehat m_K(y) = \begin{cases}1, & y \in A(Y,K),\\ 0, & y \notin A(Y,K).
\\
\end{cases}
$$

Denote by $E_x$  the degenerate distribution concentrated at a
point $x \in X$. For  $\mu \in {M^1}(X)$ denote by $\sigma(\mu)$
the support of $\mu$.

\bigskip

{\bf 2. Solution of Problem 1.}
Let $\xi_1, \ \xi_2$ be  independent  random variables with values
in the group $X$ and with distributions $\mu_1, \mu_2$. It is easily
seen that the sum  $\xi_1+\xi_2$ and the difference  $\xi_1-\xi_2$
are independent if and only if the characteristic functions
 $\widehat\mu_k(y)={\bf E}[(\xi_k, y)]$
satisfy the equation
\begin{equation}
\label{3}\widehat\mu_1(u+v)\widehat\mu_2(u-v)=
\widehat\mu_1(u)\widehat\mu_2(u)\widehat\mu_1(v)\widehat\mu_2(-v), \quad u, v \in
Y.
\end{equation}

In the sequel we need the following lemmas.

{\bf Lemma 1} (\cite{Fe1}, see also  (\cite[\S 9]{Fe4}) {\it Let
$\xi_1$ and $\xi_2$ be independent random variables with values in
a group  $X$ and
 distributions $\mu_1$ and $\mu_2$. If the sum $\xi_1+\xi_2$ and
the difference $\xi_1-\xi_2$ are independent, then distributions
$\mu_k$ can be replaced by their shifts $\mu'_k$ in such a manner
that $\sigma(\mu'_k)\subset M, \ k=1, 2$, where $M$ is a subgroup
of
 $X$ such that  $M$ is topologically isomorphic to a group
 of the form
$\mathbf{R}^m\times K$, where $m \ge 0$ and $K$ is a compact
Corwin group.}

{\bf Lemma 2}  (\cite{Fe6}). {\it Let a locally compact Abelian
group $X$ be of the form $X=\mathbf{R}^m\times K$, where $m \ge 0$
and $K$ is a compact Corwin group. Then
$\overline{Y^{(2)}}={Y^{(2)}}$ and $X_{(2)} \subset c_X$.}

{\bf Lemma 3} (\cite{Fe6}).  {\it Let $X$ be a
finite Abelian group. Then for any   function  $f(y)$ on the
group $Y$ there is a complex measure
 $\delta$ on $X$ such that $\widehat\delta (y)=f(y), \ y \in Y.$
 If  all nonzero elements of   $X$ have   order $2$, and
 $f(y)$ is a real-valued
function, then  $\delta$ is a signed measure.}

{\bf Lemma  4} (\cite{Fe5}).  {\it Let $Y$ be a locally compact
Abelain group, $\psi(y)$ be a continuous function on $Y$
satisfying the equation
$$
\Delta^2_{h}\Delta_{2k} \psi(y) = 0, \quad h, k, y \in Y,
$$ and the conditions $\psi(-y)=\psi(y), \ \psi(0)=0.$ Let $$Y = \bigcup_\alpha
{(y_\alpha + \overline{Y^{(2)}})}$$ be the
decomposition of the group $Y$ with respect to the subgroup $\overline{Y^{(2)}}$. Then
the function $\psi(y)$ can be represented in the form
$$ \psi(y) = \varphi(y) + r_\alpha, \quad y \in
y_\alpha + \overline{Y^{(2)}},
$$
\noindent where $\varphi(y)$  is a continuous function on $Y$
satisfying equation $(\ref{2})$.}

The main result of the article is the following theorem.

{\bf Theorem 1}. {\it Let $X$ be a   second countable locally
compact Abelian group. Assume that the connected component of zero
of  the group $X$ contains a finite number elements of order $2$.
Let $\xi_1$ and $\xi_2$ be independent random variables with
values in    $X$ and
 distributions $\mu_1$ and $\mu_2$ such that the sum $\xi_1+\xi_2$ and
the difference $\xi_1-\xi_2$ are independent. Then the following
statements hold.

$1$. There exists a subgroup $G \subset X_{(2)}$ such that
distributions $p(\mu_k)$ ($p$ is the natural homomorphism $p:X
\mapsto X/G$) are of the form:
$$
p(\mu_k)=\gamma*\pi_k*m_V*E_{x_k},
$$
where $\gamma \in\Gamma^s(X/G),
\ \pi_k$ are signed measures on $\left(X/G\right)_{(4)}, \ V$ is a
compact Corwin subgroup of the factor-group $X/G, \ x_k \in X/G$.

$2$. If a group $X$ is topologically isomorphic to a group of the
form  $\mathbf{R}^m\times K$, where $m \ge 0$ and $K$ is a compact
Corwin group, then either $\widehat\mu_1(y) \equiv 0$ or
$\widehat\mu_2(y) \equiv 0$  for each coset $y_j + {Y^{(2)}}$ which
disjoint with $A(Y, G)$.}

{\bf Proof}. By Lemma 1 the distributions $\mu_k$ can be replaced
with their shifts $\mu'_k$ in such a manner that
$\sigma(\mu'_k)\subset M$, where $M$ is a subgroup of
 $X$ such that  $M$ is topologically isomorphic to a group
 of the form
$\mathbf{R}^m\times K$, where $m \ge 0$ and $K$ is a compact
Corwin group.  Thus, we may assume from the beginning that the
group $X$ is of the mentioned form. Then by Lemma  2 $X_{(2)}
\subset c_X$, and hence $X_{(2)}$ is a finite subgroup. We note
that the character group $Y$ is topologically isomorphic to a
group of the form $\mathbf{R}^m\times D$, where $D$ is a discrete
group without elements of order 2. It is obvious that
$\overline{Y^{(2)}}={Y^{(2)}}$. Let $|X_{(2)}|=n$, and
\begin{equation}
\label{a1}Y=\bigcup\limits_{j=0}^{n-1}(y_j+{Y^{(2)}})
\end{equation}
be a decomposition of the group $Y$ with respect to the subgroup
${Y^{(2)}}$. Put ${N_k=\{y \in Y: \widehat\mu_k(y) \ne 0\}}, \
k=1, 2, \quad N=N_1 \cap N_2$. It follows from equation (\ref{3})
that $N$ is an open subgroup of $Y$  satisfying the following
conditions: if $2y\in N$, then $y\in N$, and
\begin{equation}
\label{a2}N\cap{Y^{(2)}}={N^{(2)}}.
\end{equation}
Equation (\ref{3}) also implies that
$$|\widehat\mu_1(u+v)| |\widehat\mu_2(u-v)|= |\widehat\mu_1(u-v)|
|\widehat\mu_2(u+v)|
$$
for all $u, v \in Y$.
  It follows from this that for arbitrary
elements $a$ and  $b$ from a given coset $y_j+{Y^{(2)}}$ we have
the equality
 \begin{equation}
\label{4}|\widehat\mu_1(a)| |\widehat\mu_2(b)|= |\widehat\mu_1(b)|
|\widehat\mu_2(a)|.
\end{equation}
Denote by  $H$ a union of   cosets $y_j+{Y^{(2)}}$ such that
$N \cap (y_j+{Y^{(2)}}) \ne \emptyset$. Since $N$ is a subgroup of
$Y$, we conclude that $H$ is also a  subgroup of $Y$. Changing if it is necessary
the numeration,
we can assume  that
\begin{equation}
\label{a3}H =\bigcup\limits_{j=0}^{l-1} (y_j+{Y^{(2)}}).
\end{equation}
Moreover, we may  assume that $y_j \in N$. Note also that
(\ref{a2}) implies that
\begin{equation}
\label{a4}N=\bigcup\limits_{j=0}^{l-1}(y_j+{N^{(2)}}).
\end{equation}

Put $G=A(X, H)$. Then $H^* \cong X/G$ and $H=A(Y, G)$. It is clear
that $G \subset X_{(2)}$. Assume that $N \cap (y_j+{Y^{(2)}}) =
\emptyset$  for some $j$. If $a \in N_1 \cap (y_j+{Y^{(2)}})$,
then
 $a \notin N_2 \cap (y_j+{Y^{(2)}})$, and the right-hand side of
(\ref{4}) vanishes. Hence, $\widehat\mu_2(b)=0$ for any $b
\in y_j+{Y^{(2)}}$. Thus, we proved  statement 2 of Theorem 1. This
reasoning  also shows that if $N \cap (y_j+{Y^{(2)}}) \ne
\emptyset$,
 then $N_1 \cap (y_j+{Y^{(2)}}) = N_2 \cap
(y_j+{Y^{(2)}})$.

 Consider the restriction of equation  (\ref{3}) to $N$.
It is obvious that the functions $|\widehat\mu_1(y)|$ and
$|\widehat\mu_2(y)|$ also satisfy equation (\ref{3}). Put $f_k(y)=
-\ln|\widehat\mu_k(y)|, \ \ y \in N, \ k=1, 2$. It follows from
(\ref{3}) that the functions  $f_k(y)$ satisfy the equation
\begin{equation}
\label{5} f_1(u+v) + f_2(u-v) = A(u) + A(v), \quad u, v \in N,
\end{equation}
where $A(u) = f_1(u) + f_2(u)$. Apply the finite difference method
to solve equation (\ref{5}). Let $h$ be an arbitrary element of
$N$. Substitute $u + h$ for $u$ and $v+h$ for $v$  in equation
  (\ref{5}). Subtracting
equation (\ref{5}) from the resulting equation we obtain
\begin{equation}
\label{6}\Delta_{2h}f_1(u+v)=\Delta_h A(u) + \Delta_h A(v), \quad
u,v
 \in N.
\end{equation}
Put in (\ref{6}) $v=0$ and subtract from (\ref{6}) the obtained
equation. We get
$$
\Delta_v\Delta_{2h}f_1(u)=\Delta_h A(v) + \Delta_h A(0),
 \quad u,v \in N.
$$
It follows from this that
$$
\Delta_v^2\Delta_{2h}f_1(u)=0, \quad u,v, h \in N.
$$
The function $f_2(y)$ satisfies the same equation. Applying Lemma
4 we obtain the representations:
$$
f_k(y)=\varphi_k(y) + p_{k, j}, \quad y \in y_j +{N^{(2)}}, \quad k=1,
2.
$$
Since equation (\ref{3}) implies that
$|\widehat\mu_1(2y)|=|\widehat\mu_2(2y)|, \ y \in Y$, the
functions $\varphi_1(y)$ and $\varphi_2(y)$ coincide on $N^{(2)}$
and hence, they also coincide on $N$, i.e.
$$\varphi_1(y)=\varphi_2(y)=\varphi(y), \ y \in N.$$
Taking this into account, we get from (\ref{3}) that
$p_{1, j}=-p_{2, j}=p_j, \ j=0, 1, \dots, l-1$. It is obvious that
$\varphi(y) \ge 0$.

 Put $l_k(y)=\widehat\mu_k(y)/|\widehat\mu_k(y)|, \ y \in N, \
k =1, 2$. The functions $l_1(y)$ and $l_2(y)$ satisfy the equation
\begin{equation}
\label{7}l_1(u+v)l_2(u-v)=l_1(u)l_2(u)l_1(v)l_2(-v), \quad u,v \in
N.
\end{equation}
Moreover, $|l_k(u)|=1, \ l_k(-u)=\overline{l_k(u)}, \ u \in N, \
l_k(0)=1, \ k=1, 2$. Putting in (\ref{7}) $v=u$, and then $v=-u$,
we find
\begin{equation}
\label{8}l_k(2u)=l_k^2(u), \quad u \in N, \quad k=1, 2.
\end{equation}

Replacing in equation (\ref{7})   $u$ by $v$ we get
\begin{equation}
\label{9}l_1(u+v)l_2(v-u)=l_1(v)l_2(v)l_1(u)l_2(-u), \quad u,v \in
N.
\end{equation}
Multiplying (\ref{7}) and (\ref{9}) we find
\begin{equation}
\label{10}l^2_1(u+v)=l^2_1(u)l^2_1(v), \quad u,v \in N.
\end{equation}
It follows from (\ref{8}) and (\ref{10}) that the restriction of
the function $l_1(y)$ to ${N^{(2)}}$ is a character of the group
${N^{(2)}}$. The same reasoning shows that  the restriction of the
function $l_2(y)$ to ${N^{(2)}}$ is also a character of the group
${N^{(2)}}$. Extend these characters from
 ${N^{(2)}}$  to some characters of the group $H$.
By the duality theorem there exist elements  $x_k \in X/G$ such
that
 $l_k(y)=(x_k, y), \ y \in {N^{(2)}}.$ Put
$l'_k(y)=(-x_k, y) l_k(y), \ y \in H$. Then
\begin{equation}
\label{c1}l'_k(2y)=1, \ y \in N, \ k=1, 2,
\end{equation}
and (\ref{8}) implies that $l'_k(y) = \pm 1, \ y \in N, \ k=1, 2$.
Hence,
\begin{equation}
\label{n2}l'_k(y)=l'_k(-y), \quad y \in N, \quad k=1, 2,
\end{equation}
and
the functions $l'_1(y)$ and $l'_2(y)$ satisfy the equation
\begin{equation}
\label{11}l'_1(u+v)l'_2(u-v)=l'_1(u)l'_2(u)l'_1(v)l'_2(v), \quad
u,v \in N.
\end{equation}
This implies that for any elements $a$ and $b$ from a given coset
$y_j + {N^{(2)}}$ the equality
\begin{equation}
\label{12}l'_1(a)l'_2(b)=l'_1(b)l'_2(a)
\end{equation}
holds. Fix an element $a \in y_j + {N^{(2)}}$. It follows from
(\ref{12}) that for all $b \in y_j + {N^{(2)}}$ either $l'_1(b)
\equiv l'_2(b)$ or
 $l'_1(b) \equiv -l'_2(b)$.
Consider the subgroup $$L={N^{(2)}}\cup(y_j+{N^{(2)}}),$$ where
 $y_j \in N$. Taking into account that $l'_1(y) \equiv l'_2(y)$ for $y \in
 {N^{(2)}}$,
 and for $y \in y_j + {N^{(2)}}$
either $l'_1(y) \equiv l'_2(y)$ or $l'_1(y) \equiv -l'_2(y)$, it
easily follows from (\ref{11}) that
\begin{equation}
\label{n1}l'_k(u+v)l'_k(u-v)=1,
\ u, v \in L, \ k=1, 2.
\end{equation}
We   conclude now from (\ref{a4}) and (\ref{n1}) that
\begin{equation}
\label{13}l'_k(u+4v)=l'_k(u), \quad u, v \in N, \quad k=1, 2,
\end{equation}
i.e. the functions $l'_k(y)$ are invariant with respect to the
subgroup ${N^{(4)}}$.

It follows from representation (\ref{a1}) that $
{Y^{(2)}}=\bigcup\limits_{j=0}^{n-1}(2y_j+{Y^{(4)}})$. Taking into
consideration (\ref{a3}) this implies that
\begin{equation}
\label{nn1}H= \bigcup\limits_{0 \le i \le n-1, \ 0 \le j \le
l-1}(2y_i+y_j +{Y^{(4)}}).
\end{equation}
Similarly, (\ref{a4}) implies that
$$ N=
\bigcup\limits_{0 \le i, j \le l-1}(2y_i+y_j +{N^{(4)}}).
$$
Taking into account (\ref{nn1}), extend the functions $l'_k(y)$
from $N$ to some functions $\widetilde l'_k(y)$ on $H$ by the formulas
$$\widetilde l'_k(2y_i+y_j+u)=l'_k(2y_i+y_j), \ u \in {Y^{(4)}}, \
0 \le i, j \le l-1, $$$$ \widetilde l'_k(2y_i+y_j+u)=1, \ u \in
{Y^{(4)}}, \ l \le i \le n-1, \ 0 \le j \le l-1, \quad k=1, 2. $$
It follows from (\ref{n2}) and (\ref{nn1}) that the functions
$\widetilde l'_k(y)$ also satisfy the condition
\begin{equation}
\label{n3}\widetilde l'_k(y)=\widetilde l'_k(-y), \quad y \in
H, \quad k=1, 2.
\end{equation}
By   construction the functions $\widetilde l'_k(y)$ are invariant
with respect to the subgroup $Y^{(4)}$, in particular they satisfy
the equation
\begin{equation}
\label{13}\widetilde l'_k(u+4v)=\widetilde l'_k(u), \quad u, v \in
H, \quad k=1, 2,
\end{equation}
and hence, they define some functions on the factor-group
$H/Y^{(4)}$. Set $ F=  A(X/G, Y^{(4)}) \subset
\left(X/G\right)_{(4)} $ and note that $F\cong(H/Y^{(4)})^*$.
Taking into account that $H/Y^{(4)}$ is a finite  group and
applying Lemma 3 we conclude that there exist   complex measures
 $\delta_k$ on $X/G$ supported  in the group
$F$ such that $\widehat\delta_k(y)=\widetilde l'_k(y), \ y \in H$.
It is obvious that $F \cong (\mathbf{Z}(4))^m$ for some  $m$. An
arbitrary character of the group  $(\mathbf{Z}(4))^m$ is of the
form
$$
((k_1, \dots, k_m), (l_1,\dots, l_m))= \exp \left\{{i \pi \over
2}\sum_{j=1}^{m}k_j l_j \right\},$$ $(k_1, \dots, k_m) \in
(\mathbf{Z}(4))^m, \quad (l_1,\dots, l_m) \in
((\mathbf{Z}(4))^m)^*$.  The  complex measures $\delta_k$,
considering as complex measures on $F$, are defined by the formulas
\begin{equation}
\label{14}\delta_k \{(k_1, \dots, k_m)\}={1 \over 4^m}
\sum_{(l_1,\dots, l_m) \in ((\mathbf{Z}(4))^m)^*} \widetilde l'_k
(l_1,\dots, l_m) \exp \left\{{i \pi \over 2}\sum_{j=1}^{m}k_j l_j
\right\}, \quad k =1, 2.
\end{equation}
Note that (\ref{n3}) and (\ref{13}) imply that
\begin{equation}
\label{16}\widetilde l'_k(3y)=\widetilde l'_k(y), \quad y \in H,
\quad k=1, 2.
\end{equation}
We also note that $\exp \left\{{i \pi \over 2}\sum_{j=1}^{m}k_j
l_j \right\}= \pm i$ if and only if when $\sum_{j=1}^{m}k_j l_j$
is an odd number. Since (\ref{16}) implies that
\begin{equation}
\label{b4}\widetilde l'_k (l_1,\dots, l_m)=\widetilde l'_k
(3l_1,\dots, 3l_m),
\end{equation}
 and the numbers
 $\exp \left\{{i \pi
\over 2}\sum_{j=1}^{m}k_j l_j \right\}$ and $\exp \left\{{i \pi
\over 2}\sum_{j=1}^{m}3k_j l_j \right\}$ are complex conjugate, it
follows from (\ref{14}) that all numbers $\delta_k \{(k_1, \dots,
k_m)\}$ are real, i.e. actually $\delta_k$ are signed measures.

Return to the representations of the functions  $f_k(y)$ on $N$.
We have
$$
f_1(y)= \varphi(u) + p_j, \quad f_2(y)= \varphi(y) - p_j, \quad y
\in y_j+{N^{(2)}}, \quad j=0, 1,\dots,l-1.
$$
By Lemma  3, there exist  signed measures  $\epsilon_1$ and
$\epsilon_2$ on $\left(X/G\right)_{(2)}$ such that
\begin{equation}
\label{z1}\widehat\epsilon_1(y) = e^{p_j}, \quad
\widehat\epsilon_2(y) = e^{-p_j}, \quad
 y \in y_j + {Y^{(2)}}, \ j = 0, 1, \dots, l-1.
\end{equation}
Put $\pi_k=\delta_k*\epsilon_k, \ k=1, 2.$ Then $\pi_k$ are signed
measures on $\left(X/G\right)_{(4)}$. Extend the function
$\varphi(y)$ from $N$ to $H$ retaining its properties (\cite[\S
5.2]{He}). Denote by $\widetilde\varphi(y)$ the extended function.
Let $\gamma$ be a symmetric Gaussian distribution on the
factor-group $X/G$ with the characteristic function
$\widehat\gamma(y) = \exp\{-\widetilde\varphi(y)\}, \ y \in H$.
Set $V=A(X/G, N)$. It is easily seen that $V$ is a compact Corwin
subgroup. Thus, we obtained that the  restriction  to $H$ of the
characteristic functions of the distributions $\mu_k$ can be
represented in the form
$$\widehat\mu_k(y) =\begin{cases}
\exp\{-\varphi(y)\} \widehat\pi_k(y)
(x_k, y), & y \in N, \\ 0, & y \in H \backslash N,\\
\end{cases}
$$
$k = 1, 2$. The desired representation
$p(\mu_k)=\gamma*\pi_k*m_V*E_{x_k}, \ k=1, 2,$ results now from the
following general remark: if $\mu \in M^1(X)$, and $H$ is a closed
subgroup of $Y$, then the restriction of the characteristic
function $\widehat\mu(y)$ to $H$ is the characteristic function of
the distribution  $p(\mu)$ on the factor-group $X/G$, where
$G=A(X, H)$, and $p$ is the natural homomorphism  $p:X \mapsto X/G$.
The theorem   is proved completely.

\bigskip

{\bf 3. Some remarks.} Give some comments to Theorem 1.

{\bf Remark 1.} Let $X$ be a    second countable locally compact
Abelian group such that $X$ is topologically isomorphic to a group
of the form  $\mathbf{R}^m\times K$, where $m \ge 0$ and $K$ is a
compact Corwin group. Assume also that $X$ contains only one
element of order 2, i.e.  $X_{(2)} \cong \mathbf{Z}(2)$. Then we
can strengthen Theorem 1.  We reason as in the proof of Theorem 1
and retain the same notation.

Since  $X_{(2)}$ is a finite group and $Y^{(2)}=A(Y, X_{(2)})$, we
have $X_{(2)}\cong(X_{(2)})^*\cong Y/{Y^{(2)}}$. Hence,  a
decomposition of the group $Y$ with respect to the subgroup
 ${Y^{(2)}}$ is of the form
$Y={Y^{(2)}}\cup(y_1+{Y^{(2)}})$, and there are two possibilities
for the subgroup $H$: either $H=Y$, and then $G=A(X, Y)=\{0\}$, or
$H=Y^{(2)}$, and then $G=A(X, Y^{(2)})=X_{(2)}$.

{\bf 1}. Let $H=Y$, then
\begin{equation}
\label{d1}N={N^{(2)}}\cup(y_1+{N^{(2)}}).
\end{equation}
Since $l'_1(y)=l'_2(y)$ for $y\in N^{(2)}$ and $l'_1(y)=\pm
l'_2(y)=\pm 1$ for $y\in y_1+N^{(2)}$, it follows from (\ref{d1})
that the functions $l'_k(y)$ are characters of the subgroup $N$.
Extend these characters from
 ${N}$  to some characters of the group $Y$. By the duality theorem
 there exist elements  $t_k \in X$ such that
 $l'_k(y)=(t_k, y), \ y \in N.$ Set $z_k=x_k+t_k$. We obtain as a result the representation
 $$\mu_k=\gamma*\epsilon_k*m_V*E_{z_k}, \quad k=1, 2,$$
where $\gamma \in\Gamma^s(X), \ \epsilon_k$ are signed measures on
$ X_{(2)}$, \ $V$ is a compact Corwin subgroup of the group $X, \
z_k \in X$. Moreover, it follows from (\ref{z1}) that  such that
$\epsilon_1*\epsilon_2=E_0.$

{\bf 2}. Let   $H=Y^{(2)}$. Then $N\subset Y^{(2)}$, and it
follows from (\ref{a2}) that
\begin{equation}
\label{d2}N=N^{(2)}.
\end{equation}
There are two possibilities for the subgroup $N$: either
$N\ne\{0\}$ or $N=\{0\}$.

\text{\boldmath $a$}. Assume that $N\ne\{0\}$. Put $W=N^*$. It follows from
(\ref{d2}) that the group $W$ contains no elements of order 2.
Consider the restriction of equation (\ref{3}) to $N$ and apply
Theorem B to the group $W$. We get that the restrictions of the
characteristic functions  $\widehat\mu_k(y)$ to $N$ are the
characteristic functions of some Gaussian distributions. It follows
easily from this that
 \begin{equation}
\label{d3}p(\mu_k)=\gamma*m_V*E_{x_k}, \quad k=1, 2,
\end{equation}
 where $\gamma
\in\Gamma^s(X/X_{(2)})$, $V$ is a compact Corwin subgroup of the
factor-group $X/X_{(2)}, \ x_k \in X/X_{(2)}$. Moreover, since
either $\widehat\mu_1(y) \equiv 0$ or $\widehat\mu_2(y) \equiv 0$
for $y_1 + {Y^{(2)}}$, it is not difficult to prove taking into
account (\ref{d3}) that at least one of the distributions  $\mu_k$
is represented in the form
$$
\mu_k=\lambda*m_U*E_{z_k},
$$
 where $\lambda
\in\Gamma^s(X)$, $U$ is a compact Corwin subgroup of $X, \ z_k \in
X$.

\text{\boldmath $b$}. Suppose that $N=\{0\}$. Obviously, it is possible only if
 $X$ is a compact group. We have $V=A(X/X_{(2)}, N)=X/X_{(2)}$.
Taking into account that the equality
 $N_1 \cap {Y^{(2)}} = N_2 \cap {Y^{(2)}}=N  \cap
{Y^{(2)}}$ always holds, we get the representation
\begin{equation}
\label{d4}p(\mu_k)=m_{{X}/{X}_{(2)}}, \ k=1, 2.
\end{equation}
 Moreover, since either $\widehat\mu_1(y) \equiv 0$
or $\widehat\mu_2(y) \equiv 0$ for $y_1 + {Y^{(2)}}$, it follows
from (\ref{d4}) that $\mu_k=m_X$ at least for one of the
distributions  $\mu_k$. Let for example, $\mu_1=m_X$. Returning to
the random variables $\xi_k$, it means that $\xi_1$ and $2\xi_2$
are identically distributed random variables with distribution
$m_{X}$.

As a result we obtain the following statement.

{\bf Theorem 2}. {\it Let $X$ be a    second countable locally
compact Abelian group such that $X$ is topologically isomorphic to
a group of the form  $\mathbf{R}^m\times K$, where $m \ge 0$ and
$K$ is a compact Corwin group. Assume also that the group $X$ contains only
one element of order $2$.  Let $p$ be the natural homomorphism
$p:X \mapsto X/X_{(2)},$ and $Y={Y^{(2)}}\cup(y_1+{Y^{(2)}})$ be
a decomposition of the group $Y$ with respect to the subgroup
$Y^{(2)}$. Assume that $\xi_1, \ \xi_2$ are independent random
variables with values in    $X$ and
 distributions $\mu_1, \mu_2$ such that the sum $\xi_1+\xi_2$ and
the difference $\xi_1-\xi_2$ are independent. Then either
$$\mu_k=\gamma*\epsilon_k*m_V*E_{z_k},$$ where $\gamma \in\Gamma^s(X), \
\epsilon_k$ are signed measures on $ X_{(2)}$ such that
$\epsilon_1*\epsilon_2=E_0,$ \ $V$ is a compact Corwin subgroup of
$X, \ z_k \in X, \  k=1, 2,$ or
$$p(\mu_k)=\gamma*m_V*E_{x_k},$$
 where $\gamma
\in\Gamma^s(X/X_{(2)})$ $V$ is a compact Corwin subgroup of  the
factor-group $X/X_{(2)}, \ x_k \in X/X_{(2)}, \  k=1, 2$, and at
least one of the distributions  $\mu_k$ is represented in the form
$$
\mu_k=\lambda*m_U*E_{z_k},
$$
 where $\lambda
\in\Gamma^s(X)$, $U$ is a compact Corwin subgroup of $X, \ z_k \in
X$.}

We illustrate Theorem 2 by the following example. Let
$X=\mathbf{T}$. Then $Y \cong \mathbf{Z}$. In order not to complicate notation
we assume that $Y=\mathbf{Z}$. We have
$\mathbf{T}_{(2)}\cong\mathbf{Z}(2)$, and hence there are two possibilities:
either $G=\{0\}$ or $G=\mathbf{T}_{(2)}$.

{\bf 1}. $G=\{0\}$.  We observe that all compact Corwin subgroups
$V$ of the group $\mathbf{T}$ are of the form: either
$V=\mathbf{T}$, or $V$ is the subgroup of $m$th roots of 1, where
$m$ is odd, i.e. $V\cong \mathbf{Z}(m)$. If $V=\mathbf{T}$, then
by Theorem $2$ $$\mu_1=\mu_2=m_\mathbf{T}.$$ Let $V \cong
\mathbf{Z}(m)$. Then  by Theorem 2
\begin{equation}
\label{b1}\mu_k=\gamma*\epsilon_k*m_V*E_{z_k},
\end{equation}
 where $\gamma
\in\Gamma^s(\mathbf{T}), \ \epsilon_k$ are signed measures on
 $ \mathbf{T}_{(2)}$ such that $\epsilon_1*\epsilon_2=E_1,$   $ \ z_k \in
\mathbf{T}, \  k=1, 2.$  It follows from $V \cong \mathbf{Z}(m)$  that
$N=\mathbf{Z}^{(m)}$, and hence (\ref{b1}) implies that
the characteristic functions of distributions
 $\mu_k$ are represented in the form
$$
\widehat\mu_1(n) = \begin{cases} \exp\{-\sigma n^2 + int_1 \},
& n \in \mathbf{Z}^{(2m)}, \\ \exp\{-\sigma n^2
+ int_1 + q \}, & n \in \mathbf{Z}^{(m)}\backslash
\mathbf{Z}^{(2m)},\\  0,
& n \notin
\mathbf{Z}^{(m)}, \\
\end{cases} $$

$$
\widehat\mu_2(n) = \begin{cases} \exp\{-\sigma n^2 + int_2 \},
& n \in \mathbf{Z}^{(2m)}, \\ \exp\{-\sigma n^2
+ int_2 - q \}, & n \in \mathbf{Z}^{(m)}\backslash \mathbf{Z}^{(2m)},
\\ 0, &
n \notin \mathbf{Z}^{(m)}, \\
\end{cases} $$

\noindent where $ \sigma \ge 0, \ t_k, \ q \in \mathbf{R}, \ k= 1,
2$.

{\bf 2}. $G=\mathbf{T}_{(2)}$. It follows from (\ref{d2}) that
$N=\{0\}$. Hence, $V=A(\mathbf{T}/\mathbf{T}_{(2)}, N)=
\mathbf{T}/\mathbf{T}_{(2)}$. We get by Theorem 2
$ p(\mu_k)=m_{\mathbf{T}/\mathbf{T}_{(2)}}, \ k=1, 2$, i.e.
$\widehat\mu_1(2n)=\widehat\mu_2(2n)=0, \ n \in \mathbf{Z}, \ n
\ne 0.$ Moreover, either
 $\widehat\mu_1(2n+1)= 0$ or $\widehat\mu_2(2n+1)=0, \ n \in \mathbf{Z}.$
This implies that at least one of the distributions $\mu_k$, say
$\mu_1=m_\mathbf{T}$. As regards to the second distribution, we
know only that
 $\widehat\mu_2(2n)=0, \ n \in \mathbf{Z}, \ n
\ne 0$. Returning to the random variables  $\xi_k$ it means that
 $\xi_1$ and  $2\xi_2$ are identically distributed random variables
 with the distribution
 $m_\mathbf{T}$.

The obtained description of possible distributions  $\mu_k$ for
the group $X=\mathbf{T}$ is the main content of article
\cite{BaESta}. Similarly, it is possible to get
from Theorem 2 the description of possible distributions
 $\mu_k$ for the group $\mathbf{R}\times
\mathbf{T}$ and for the \text{\boldmath $a$}-adic solenoids $\Sigma_\text{\boldmath $a$}$
found in \cite{My}.

{\bf Remark  2.} Let $X=\mathbf{T}^2$. Then
$Y\cong\mathbf{Z}^2$. Denote by $x=(e^{it}, e^{is})$ elements of the group
 $X$, and by $y=(m, n)\in
\mathbf{Z}^2$  elements of the group
 $Y$. We will construct independent random variables $\xi_1$
and $\xi_2$ with values in the group $X$ and distributions
$\mu_1$ and $\mu_2$ such that the sum $\xi_1+\xi_2$ and the
difference $\xi_1-\xi_2$ are independent, and
$\mu_k=\gamma*\pi_k$, where $\gamma\in \Gamma^s(X)$, and $\pi_k$
are
 signed measures supported in $X_{(4)}$, but not in $X_{(2)}$.
Thus, we will show that the statement in Theorem 1 that  $\pi_k$
are signed measures on $\left(X/G\right)_{(4)}$, generally
speaking can not be strengthen  to the statement that $\pi_k$ are
signed measures on
 $\left(X/G\right)_{(2)}$ (compare with Theorem  $2$).

Take $\sigma > 0$ such that
\begin{equation}
\label{b2}\sum\limits_{(m, n)\in \mathbf{Z}^2, \ (m, n)\ne (0, 0)}
\quad e^{-\sigma(m^2+n^2)} < 1.
\end{equation}
Consider on the group $Y$ the functions

$$l_1(m,n)=\begin{cases} 1, \ \ \ (m,n) \in \{Y^{(2)}, \ (1,0) +
Y^{(4)}, \ (3,0) + Y^{(4)}, \ (0,1) + Y^{(4)}, \\
\ \ \ \ \ \ (0,3) + Y^{(4)}, \ (1,1) +
Y^{(4)}, \ (3,3) + Y^{(4)} \}, \\ -1, \   (m,n) \in \{(1,2)
+ Y^{(4)}, \ (3,2) + Y^{(4)}, \ (2,1) + Y^{(4)}, \\
\ \ \ \ \ \ (2,3) + Y^{(4)}, \ (1,3) +
Y^{(4)}, \ (3,1 ) + Y^{(4)}  \}\\
\end{cases}
$$
and
$$
l_2(m,n)=\begin{cases} 1, \ \ \ (m,n) \in \{Y^{(2)}, \ (1,0) +
Y^{(4)}, \ (3,0) + Y^{(4)}, \ (0,1) + Y^{(4)}, \\
\ \ \ \ \ \ (0,3) + Y^{(4)}, \ (1,3) +
Y^{(4)}, \ (3,1) + Y^{(4)} \}, \\ -1, \
 (m,n) \in \{(1,2)
+ Y^{(4)}, \ (3,2) + Y^{(4)}, \ (2,1) + Y^{(4)}, \\
\ \ \ \ \ \ (2,3) + Y^{(4)}, \ (1,1) +
Y^{(4)}, \ (3,3 ) + Y^{(4)}  \}.\\
\end{cases}
$$

 Consider on the group $X$ the functions
\begin{equation}
\label{b3}\rho_k(x)=\rho_k(e^{it}, e^{is})=1+ \sum\limits_{(m,
n)\in \mathbf{Z}^2, \ (m, n)\ne (0, 0)} \quad
e^{-\sigma(m^2+n^2)}l_k(m,n)e^{-i(tm+sn)}, \quad k=1, 2.
\end{equation}
Inequality (\ref{b2}) implies that $\rho_k(x) > 0$ and it is
obvious that
$$
\int_X\rho_k(x)dm_X(x)=1, \quad k=1, 2.
$$
 Let $\mu_k$ be the distributions on the group  $X$ with the densities
  $\rho_k(x)$ with respect to the Haar distribution
 $m_X$, and let $\xi_k$ be independent random variables
 with values in    $X$ and distributions  $\mu_k, \ k=1, 2$.
It follows from (\ref{b3}) that
\begin{equation}
\label{z2}\widehat\mu_k(m,n)=e^{-\sigma(m^2+n^2)}l_k(m,n), \quad
(m, n)\in \mathbf{Z}^2, \quad k=1, 2.
\end{equation}
 We can check directly that the
functions
 $l_k(m, n)$ satisfy equation
 (\ref{3}), and hence the characteristic
functions $\widehat\mu_k(m,n)$ also satisfy equation (\ref{3}).
Thus, the sum $\xi_1+\xi_2$ and the difference $\xi_1-\xi_2$ are
independent.

By the construction, the functions $l_k(m, n)$ are invariant with respect to the subgroup
 $Y^{(4)}$, and hence, they define some functions
 $\widetilde l_k(m, n)$ on the factor-group  $Y/Y^{(4)}\cong (\mathbf{Z}(4))^2$.
  By Lemma 3 there exist complex measures $\pi_k$ on $X$ supported
  in $A(X, Y^{(4)})=X_{(4)}$ such that
  \begin{equation}
\label{z3}\widehat\pi_k(m,n)=  l_k(m,n), \ (m, n)\in \mathbf{Z}^2,
\quad k=1, 2.
\end{equation}
  Since the functions  $\widetilde l_k(m, n)$ satisfy the condition
 (\ref{b4}), reasoning as in the proof of Theorem 1 we get that
 the complex measures  $\pi_k$ actually are signed measures.
 Since the characteristic functions   $\widehat\pi_k(m, n)$ are not
 invariant with respect to the subgroup  $Y^{(2)}$, the supports of the  signed measures
  $\pi_k$ do not contain in $X_{(2)}$. Let $\gamma$ be a symmetric Gaussian distribution
  on $X$ with the characteristic function
  \begin{equation}
\label{z4}\widehat\gamma(m, n)=e^{-\sigma(m^2+n^2)}, \  (m, n)\in
\mathbf{Z}^2.
\end{equation}
 It follows from (\ref{z2})--(\ref{z4}) that $\mu_k=\gamma*\pi_k, \ k=1, 2.$

{\bf Remark 3.} Assume that under conditions of Theorem 1 the independent
random variables $\xi_1$ and $\xi_2$ are identically distributed, i.e. $\mu_1=\mu_2=\mu$. As has been proved in \cite{Fe6}, this implies that
$$\mu=\gamma*\pi*m_V*E_x,$$ where $\gamma \in \Gamma^s(X)$, $\pi$ is a signed measure on $X_{(2)}$ such that $\pi^{*2}=E_0$, $V$ is a compact Corwin subgroup of the group $X$, $x\in X$. Taking into account Theorem 1 and the example constructing in Remark 2, we see that even we assume that in Theorem 1 $G=\{0\}$, the distributions $\mu_k$ have generally speaking factorisation much more complicated than in the case of identically distributed random variables.

\newpage

 G.M. Feldman

\bigskip

B. Verkin Institute for Low Temperature Physics and Engineering

of the National Academy of Sciences of Ukraine

\bigskip

Mailing address:

\bigskip

B. Verkin Institute for Low Temperature Physics and Engineering

Lenin Ave. 47

61103 Kharkov, Ukraine

\bigskip

e-mail: feldman@ilt.kharkov.ua

\end{document}